# A Novel Geometric Model of Natural Spirals Based on Right Triangle Polygonal Modular Formations

Robert E. Grant[1], Talal Ghannam[2], Amanda Kennedy[3]

[1]Strathspey Crown Holdings, Crown Sterling. Newport Beach, California, USA. [2,3]Crown Sterling, Newport Beach, California, USA.

We propose a novel class of spirals that are based on perfect polygonal formations. The spirals are defined by a fractal of right triangles that delineate their geometry and determine their progression rates and modes. We show how these spirals match naturally occurring ones, including hurricanes and galaxies, much better than already proposed models, allowing for more accurate categorization of these phenomena and better prediction of their physical properties.

## I.  INTRODUCTION

The universe is defined by various patterns and forms, with spirals being the dominant among them, observed in the faraway galaxies as well as in hurricanes and tornados, plants and trees, water flow, DNA, etc. They are the archetypical patterns upon which most types of matter, living and nonliving, organize themselves. Therefore, understanding the physical properties of these fundamental shapes and the mechanisms by which they form is essential to our understanding of the whole universe, on the macroscopic level as well as on the microscopic one.

The first approach scientists usually take to understand a phenomenon is by formulating it via an appropriate mathematical model. The general research on the topic has concentrated mainly on the hyperbolic logarithmic spiral (HLS)[1], where the Euler number $e$ (2.718…) works as the expansion factor. Nevertheless, this approach does not produce accurate models that can fit the observed forms of all spirals unless it is adjusted by certain factors, depending on the natural spiral being modeled.[2,3] For example, the equation $ln(r - r_0) = A_1 + B_1(\theta)$ was proposed as a good fit for some of the rain bands of hurricanes.[4] As obvious, many parameters must be included in order to force the logarithmic spiral to match the natural one. And still, the resemblance varies from good to not that close.

Accurately categorizing natural spirals based on specific mathematical and geometrical models offers many benefits. For example, there has been an awareness of the relationship between spirals shapes and their physical properties for over five decades.[5] Thus, we can predict the intensity of a specific hurricane simply by matching its radar images to a known model that is already associated with a specific category.

Nature is intrinsically simple, and the triangle is the simplest of all geometrical shapes and also the most stable among them. Moreover, any geometric form can be deconstructed into a matrix of triangles of different shapes and sizes, as in finite element analysis. The right triangle is of particular significance due to the Pythagorean relationship that governs its sides, with the sum of the square of the two perpendicular sides equal to the square of the third one, the hypotenuse. Therefore, the right triangle is the most suitable shape to be used in constructing spirals, mathematically as well as geometrically.

Using right triangles, we were able to create spiraling forms that fit most natural spirals much better than any other model proposed in the literature. This we achieved without needing an ample parameter space that requires adjustment depending on the modeled spiral. Moreover, we will show that there is a direct correlation between the number of triangles forming one complete round of the spiral (which we call the *mod*) and the specific properties of the natural spirals that match it.

Below are the definitions of some of the terms used throughout the paper:

Rain bands: The spiral-like structure of clouds in which heavy rainfall resides.

Regular polygon: a shape with the exact same side lengths, not to be confused with a regular polygon.

Mod: The number of right triangles formed from intersecting a perfect polygon at all its vertices and midpoints.

## II.    Triangular Spiral: Math and Geometry

Forming spirals from right triangles is not a new concept. The spiral of Theodorus is one such example where right triangles, all having base lengths of 1 unit, can be stacked one next to the other to create a spiraling behavior, as shown below.

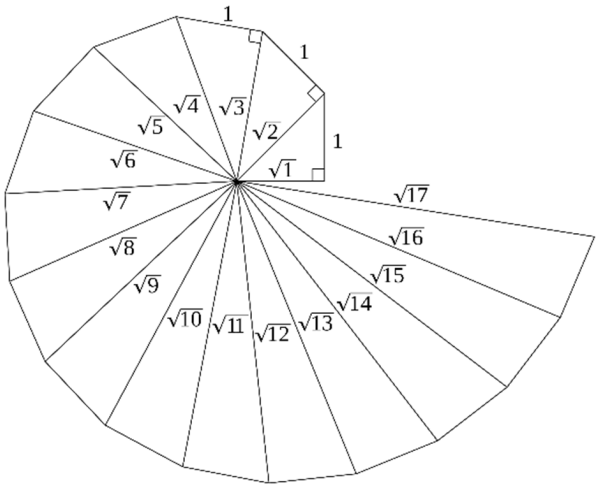

Figure 1: The spiral of Theodorus, made of a matrix of right triangles with the outer base of all triangles having a length of 1 unit.

Through our research, we were able to build other types of spirals by juxtaposing right triangles such that each triangle is constrained by the base of the previous one. The emerging shapes delineate perfect regular polygonal geometry, such as squares, pentagons, hexagons, etc.

Consider the spiral below, which requires 12 triangles for each round. Note how the base of each triangle constrains the dimension of the next by defining its height.

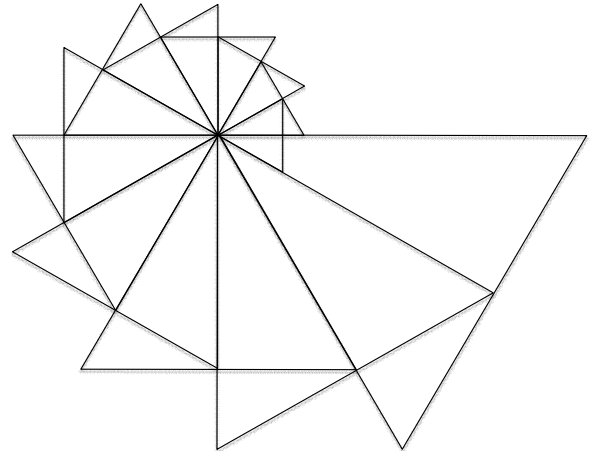

Figure 2: A new form of spirals made from right triangles juxtaposed such that the base of the previous triangle work as the height of the next one.

The triangles that create these *triangular spirals* have similar proportions, varying in scale only. For the above specific example, the dimensions have ratios of $[1:\sqrt{3}:2]$, as shown below. Thus, if we start from a triangle of $[1, \sqrt{3}, 2]$ units, the next one will have a hypotenuse of $4/\sqrt{3}$, and so on.

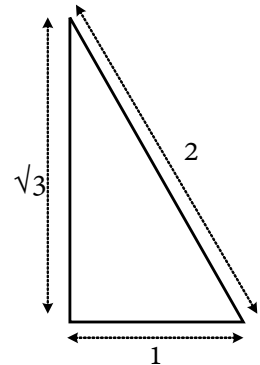

Figure 3: The right triangle that forms the base of the 12-turns spiral has the dimensions of $[1:2:\sqrt{3}]$.

The value √(4/3) = 1.1547... is the ratio of the hypotenuse to the height, which governs how each new triangle grows from the previous one. It acts as a log base (similar to *e*) and accounts for the expansion rate around a circle. The growth rate of 1.1547… is equal to 1/cos(30°) = 1.1547…, and for this specific angle, the spiral is formed based on a hexagonal structure, as shown below. Generally, the number of twists or turns (the mod) a triangular spiral performs is twice the number of sides of the polygon it delineates. So for a specific polygonal geometry of *n*-sides, the mod is 2*n*. Therefore, the above-mentioned spiral is mod(12).

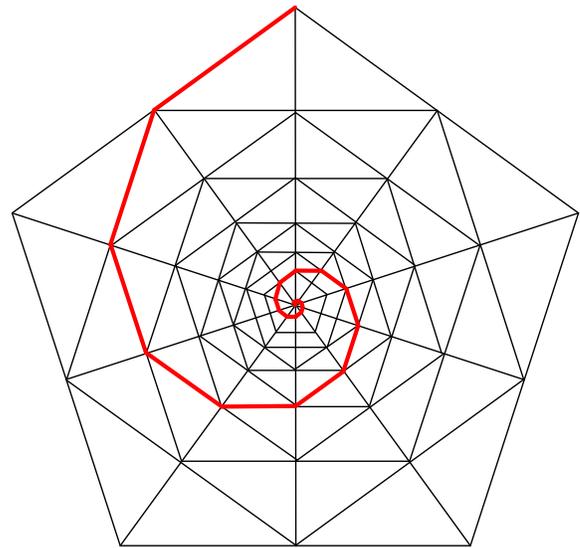

Figure 5: A pentagonal spiral emerging from connecting the tip of each right triangle to its predecessor.

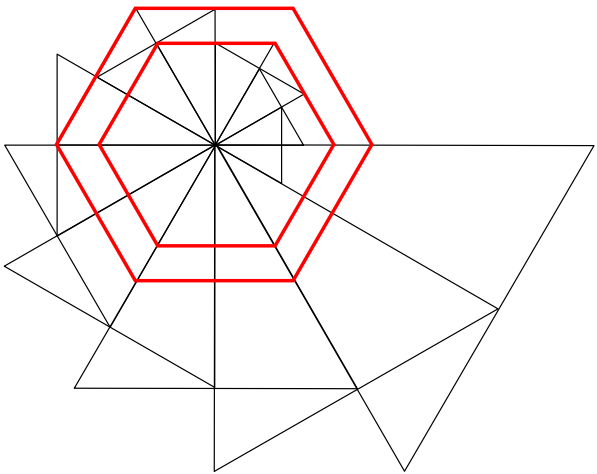

Figure 4: The spiral defined by a hexagon performs 12 twists for each 360° complete turn; thus, it is a mod(12) triangular spiral.

The higher the expansion base values, the lower the number of right triangles per loop (lower mod number), and consequently, the less the circularity of the spiral.

Another spiral emerges from a pentagonal structure, as shown below. The triangles that form this spiral have proportions [1:1.37639:1.7013], and therefore, each triangle grows from the previous one by the ratio of 1.7013…/1.37639… = 1/cos(36°) = 1.23606… The number of twists it takes to make one full round is 2×5 = 10, and this is a mod(10) spiral.

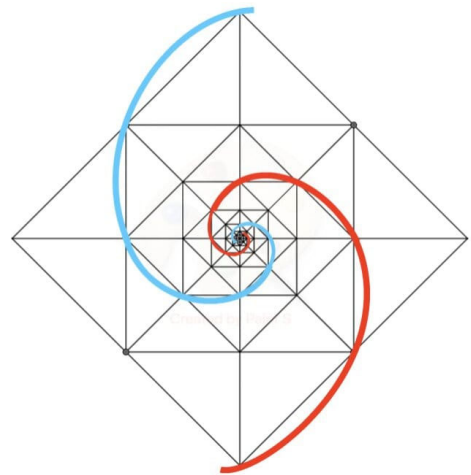

Figure 6: The square, a spiral of mod(8), displays a spiral pitch with an expansion value of 1.414213562. The spiral is visibly less circular compared to higher mod ones.

A chart containing measurements of the logarithmic base values for different spirals corresponding to polygons having sides' lengths of 3 up to 109 can be found in the Appendix. (The GeoGebra online application was used to generate each spiral to 17 digits accuracy.)

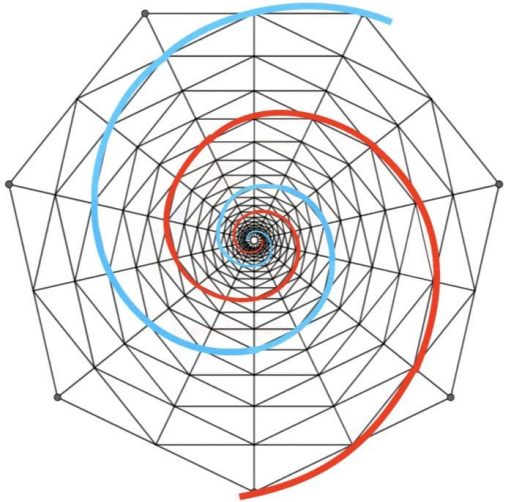

Figure 7: A heptagon spiral, mod(14), displays a spiral with a base value of 1.109916264. The higher the mod number, the more circular the spiral becomes.

Mathematically, spirals are best described through spherical or cylindrical coordinates. Each point on the spiral is determined by its distance $r$ from the origin as a function of the angle $\varphi$, or the ranking of the triangle $n$, as shown below. In our case, the radial distances coincide with the points where the extended bases intersect the hypotenuses, e.g., $a$, $r_1$, $r_2$, etc.

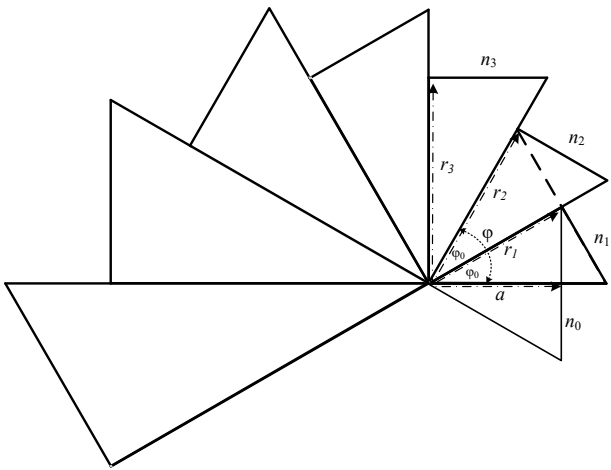

Figure 8: The triangular spiral is defined by $n$, the number of triangles, and $\varphi_0$, the angle that defines each right triangle.

As the triangles of the above spirals are all similar, differing only in scale, the internal angle $\varphi_0$ is always the same for any given spiral, differing depending on the triangles used (which determines the mod) and is restricted by the condition $\varphi_0 < 90°$ (as we can't have two 90° in one triangle). This is different from the Theodorus spiral, for example, as even though the latter is also made from right tringles, they are all different from each other, with different internal angles.

Thus, we can describe the radial coordinate of a triangular spiral as:

$$r(n) = a/\cos(\varphi_0)^n$$

Where $n$ is the order of the triangle in the sequence, having the discrete values [0, 1, 2, 3, …], and $a$ is the height of the first right triangle, the spiral seed. So, for $n = 0$, $r$ is equal to $a$. Of course, $\varphi_n$, the full angle of $r_n$ is related to $\varphi_0$ through $\varphi_n = n \cdot \varphi_0$; therefore, we can rewrite the above relationship as:

$$r(\varphi) = a/\cos(\varphi_0)^{(\varphi/\varphi_0)}$$

From here, we can find a formula for $\varphi$ by taking the natural log of both sides, which yields:

$$\varphi = k \times ln(a/r)$$

Where $k$ is a constant that depends on the triangles used in making up the spiral, defined as $\varphi_0/ln(\cos(\varphi_0))$.

The curvature of the spiral is defined by the *polar slope angle*, which measures the angle that a tangent line to the spiral at a certain point makes with the tangent line to a circle centered at the spiral's center, at that same point, as shown below.

Mathematically, the polar angle ($\alpha$) is calculated using the formula: $\tan(\alpha) = r'/r$, where $r'$ is the derivative of $r$ with respect to it argument. So, for $r(n) = a/\cos(\varphi_0)^n$, the polar angle will be:

$$\tan(\alpha) = [d/dn(a/\cos(\varphi_0)^n)]/[a/\cos(\varphi_0)^n] \to$$

$$\tan(\alpha) = -ln[\cos(\varphi_0)]$$

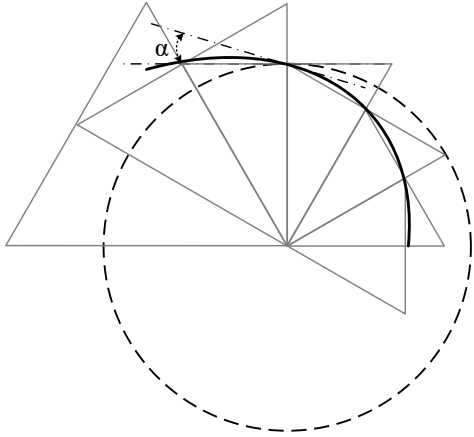

Figure 9: The polar slope angle is defined as the angle between the tangent of the spiral at a certain point and the circle that runs through the same point and is centered at the spiral origin.

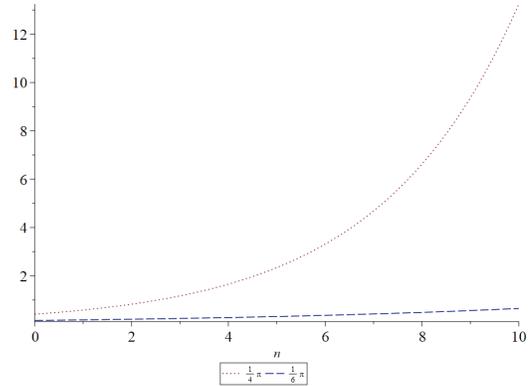

Figure 11: The radial difference $\Delta r$ between successive triangles as a function of $n$ ($a = 1$), and for two different angles $\varphi_0 = \pi/6$ (dashed line) and $\varphi_0 = \pi/4$ (dotted line).

As illustrated by the plot below, the polar slope angle diverges as we approach the critical value of $\varphi_0 = \pi/2$.

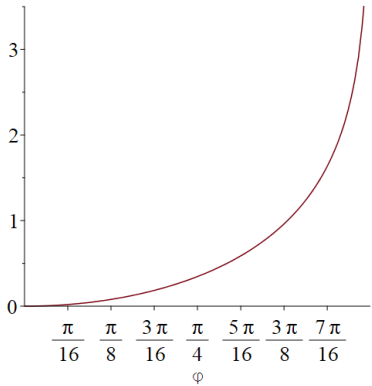

Figure 10: Polar slope angle vs. the internal angle $\varphi_0$. The slope diverges as the angle approaches $\pi/2$.

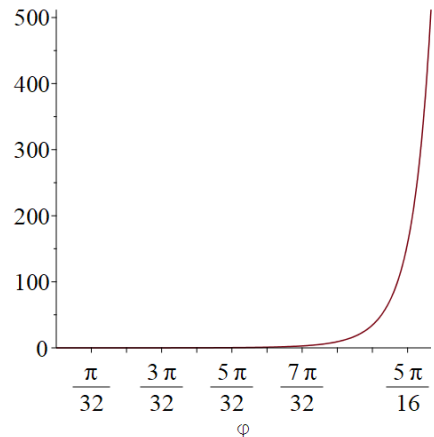

Figure 12: The radial difference diverges at around $\varphi_0 = 7\pi/32 \approx 40°$ (this is for $n = 10$).

The growth rate of the the radius of the spiral between points $n$ and $n+1$ is:

$$\Delta r = a/\cos(\varphi_0)^{n+1} - a/\cos(\varphi_0)^n \rightarrow$$

$$\Delta r = a[(1-\cos(\varphi_0))/\cos(\varphi_0)^{n+1}].$$

As shown below, for $\varphi_0 = \pi/6$, $\Delta r$ increase in an almost linear fashion with $n$. however, for larger angles, the rate of increase becomes nonlinear, reaching exponential values for values of $\varphi_0$ around 40°.

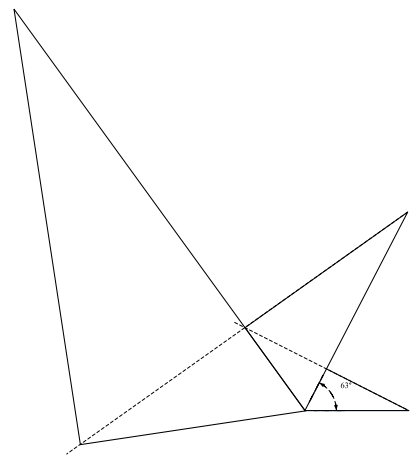

Figure 13: For large values of $\varphi_0$, the radial difference diverges quickly, as shown above for triangles having $\varphi_0 = 63°$ and for $n = [0, 1, 2]$.

## III. Triangular Spirals in Nature

We overlay triangular spirals on top of images released by NASA, National Geographic, and other scientific and educational organizations. The images are for galaxies, hurricanes, shells, etc. Apart from the opacity, no further editing or distorting of the images has been made.

The overlaying reveals an excellent match with the natural formations. For example, the nautilus shell has always been considered as a perfect example of the golden spiral working in nature (the golden spiral is defined by the irrational golden mean $\Phi = 1.618\ldots$). Nevertheless, the triangular spiral of mod(16) seems to fit the shell more perfectly than the golden spiral does, as shown below. Note the perfect match with the first three rounds of the shell. The shell then slightly deviates from the triangular spiral, but we believe as the nautilus gets older, its shell will again achieve a perfect matching. And not only does the nautilus' spiral match that of the mod(16), but also there are 16 chambers inside the shell for each turn, matching up exactly with the triangles' hypotenuses.

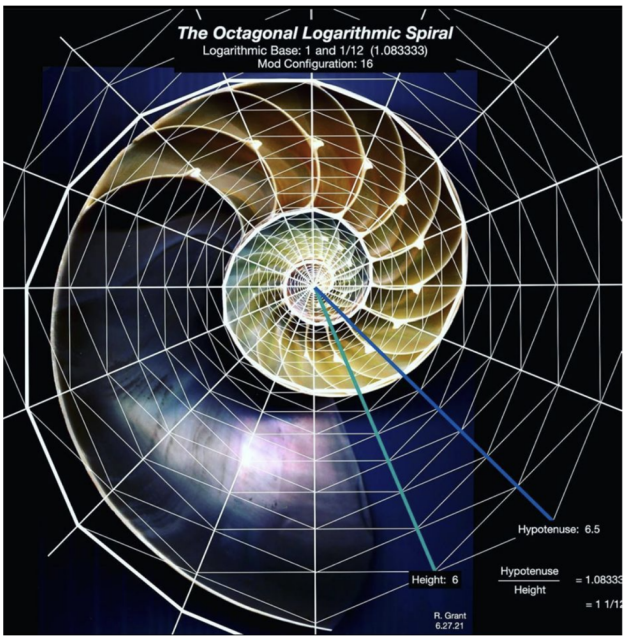

Figure 14: Overlaying a triangular spiral of mod(16) over a nautilus shell reveals a perfect matching (better than a golden spiral would do).

In the figure below, we compare a mod(8) spiral that is based on square polygonal geometry to the hurricane Sally of 2020, having a category of two. Note the perfect matching between the two spirals.

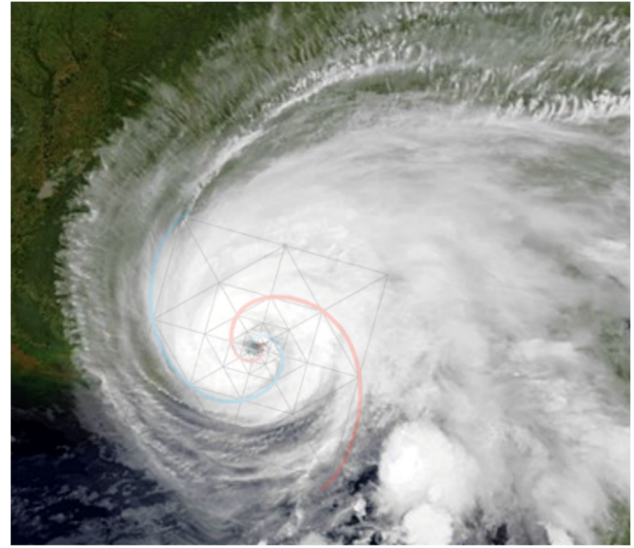

Figure 15: Hurricane Sally (category 2), whose maximum wind speed reached about 105 mph, is identified as a mod(8) spiral.

Another example is between hurricane Odile (2014) of category four and the mod(16) spiral.

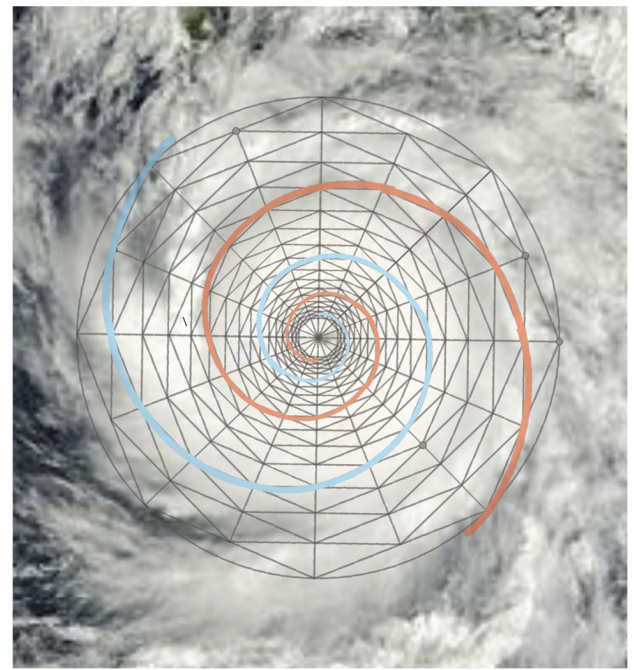

Figure 16: Hurricane Odile, category 4, matching perfectly with a triangular spiral of mod(16).

Hurricane Katrina of 2005 was a category five. It matches perfectly with a mod(18) spiral.

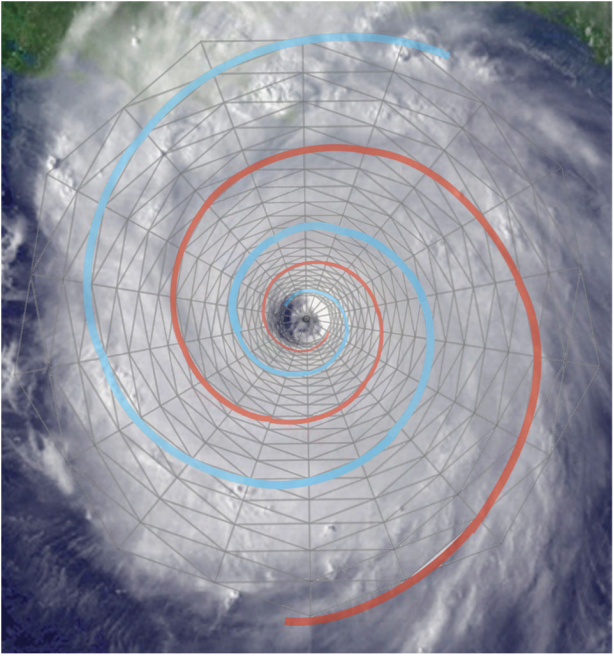

Figure 17: Hurricane Katrina. Category 5, matching a spiral of mod(18).

As demonstrated by these images, there seems to be a direct correlation between the category or the strength of the hurricane and the mod of the matching triangular spiral; the stronger the hurricane, the higher the mod. Higher mod numbers indicate a higher number of right triangles involved in creating the spiral, which is reflected by the increase in the number of sides of its respective overall polygonal geometry.

The outer space is filled with spiral-like structures as well, mainly galaxies. Galactic spiral arms are made up of different materials, including space dust, gas, celestial bodies such as stars, planets, etc. The gravitational interaction between the stars in the spiral arm, along with the effect of the dark matter that surrounds them, controls their angular momentum and hence their spirality.[6]

Galactic spirals have the biggest share of spiral-model analysis. The images below illustrate many types of spiral galaxies overlayed with their matching triangular spirals, along with their respective mods.

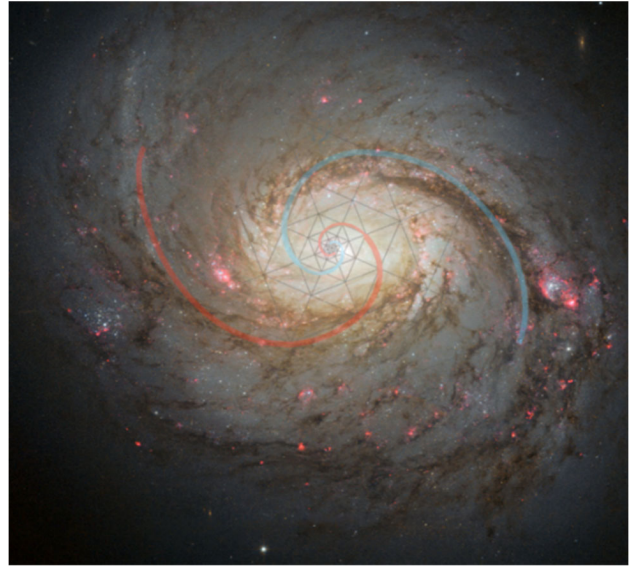

Figure 18: Galaxy NGC 1566 identified as a mod(14) spiral.

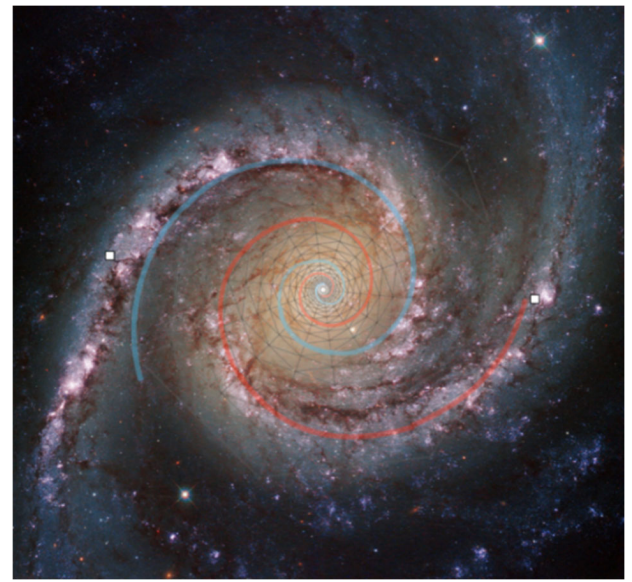

Figure 19: Galaxy M77 identified as a mod(8) spiral.

Galaxies with higher central mass concentrations are predicted to have more tightly wound spiral arms.[7] Furthermore, the greater the speed and the heavier the density, the more mods there are in the matching triangular spiral. The more mods, the tighter the spirals would wind, which creates a more circular shape.

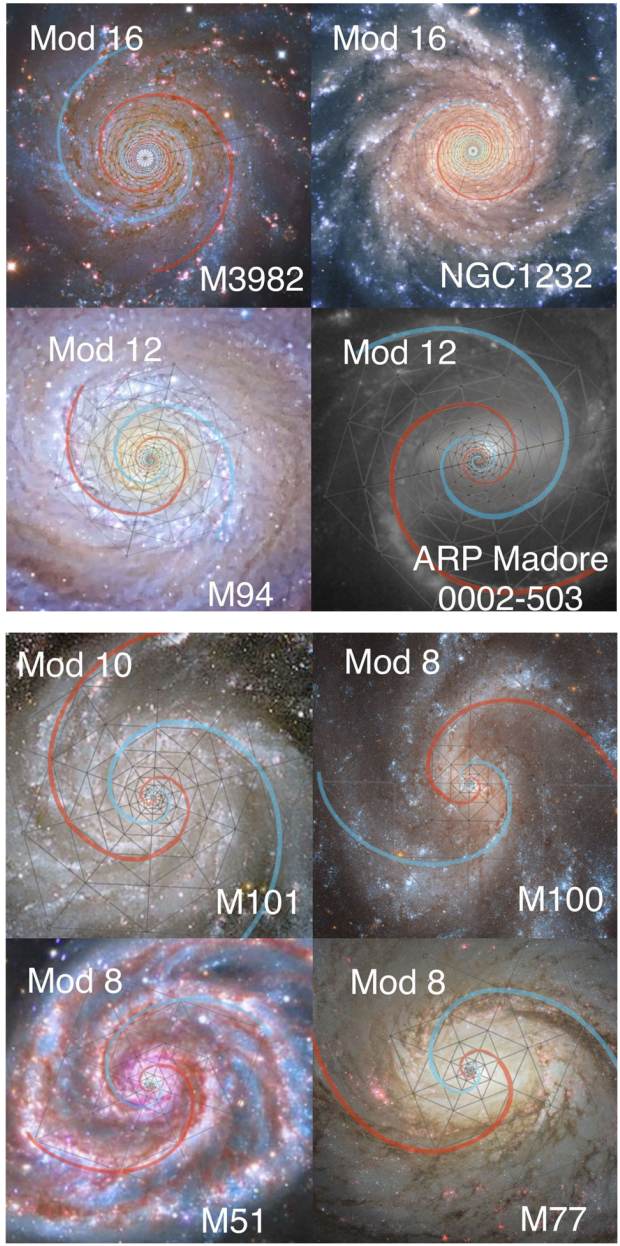

Figure 20: Another example of perfect matching between celestial galaxies and triangular spirals of different mods.

This perfect correlation allows for better categorization and prediction of the properties of the numerous spiraling formations manifested in the natural world. It also hints at a possible fundamental role the right triangular geometry plays in forming the physical universe, which is the subject of ongoing research by the authors.

## IV. Conclusion

A novel class of triangular spirals is proposed, and their geometrical and mathematical properties are analyzed. These new forms are shown to match many natural spiraling phenomena much better than any previously proposed model. Even though this study does not explain how the natural spirals are originally formed; however, distinct connections have been observed, where physical attributes of the natural spirals can be deduced based on the mods of their matching triangular spirals. An in-depth study is currently being conducted to gain a deeper insight into right triangles and how their geometry could help us better understand other fundamental physical phenomena.

**Appendix**

A chart containing all necessary measurements of the progression base values for shapes with size three up until 109.

| sides: | 0 | 1 | 2 | 3 | 4 | 5 | 6 | 7 | 8 | 9 |
|---|---|---|---|---|---|---|---|---|---|---|
| degrees: | — | — | — | 60 | 45 | 36 | 30 | 25.71428571 | 22.5 | 20 |
| 1/cos: | — | — | 2 | 1.414213562 | 1.236067977 | 1.154700538 | 1.109916264 | 1.0823922 | 1.064177772 | |

| sides: | 10 | 11 | 12 | 13 | 14 | 15 | 16 | 17 | 18 | 19 |
|---|---|---|---|---|---|---|---|---|---|---|
| degrees: | 18 | 16.36363636 | 15 | 13.84615385 | 12.85714286 | 12 | 11.25 | 10.58823529 | 10 | 9.473684211 |
| 1/cos: | 1.051462224 | 1.042217116 | 1.03527618 | 1.029927831 | 1.025716863 | 1.022340595 | 1.019591158 | 1.017321838 | 1.015426612 | 1.013827283 |

| sides: | 20 | 21 | 22 | 23 | 24 | 25 | 26 | 27 | 28 | 29 |
|---|---|---|---|---|---|---|---|---|---|---|
| degrees: | 9 | 8.571428571 | 8.181818182 | 7.826086957 | 7.5 | 7.2 | 6.923076923 | 6.666666667 | 6.428571429 | 6.206896552 |
| 1/cos: | 1.012465126 | 1.011295333 | 1.010283227 | 1.009401621 | 1.008628961 | 1.007947971 | 1.007344677 | 1.006807673 | 1.006327577 | 1.005896609 |

| sides: | 30 | 31 | 32 | 33 | 34 | 35 | 36 | 37 | 38 | 39 |
|---|---|---|---|---|---|---|---|---|---|---|
| degrees: | 6 | 5.806451613 | 5.625 | 5.454545455 | 5.294117647 | 5.142857143 | 5 | 4.864864865 | 4.736842105 | 4.615384615 |
| 1/cos: | 1.00550828 | 1.005157136 | 1.004838572 | 1.004548674 | 1.004284099 | 1.004041978 | 1.003819838 | 1.003615536 | 1.003427213 | 1.003253241 |

| sides: | 40 | 41 | 42 | 43 | 44 | 45 | 46 | 47 | 48 | 49 |
|---|---|---|---|---|---|---|---|---|---|---|
| degrees: | 4.5 | 4.390243902 | 4.285714286 | 4.186046512 | 4.090909091 | 4 | 3.913043478 | 3.829787234 | 3.75 | 3.673469388 |
| 1/cos: | 1.003092198 | 1.002942834 | 1.002804043 | 1.002674852 | 1.002554394 | 1.002441898 | 1.002336678 | 1.002238119 | 1.002145671 | 1.002058837 |

| sides: | 50 | 51 | 52 | 53 | 54 | 55 | 56 | 57 | 58 | 59 |
|---|---|---|---|---|---|---|---|---|---|---|
| degrees: | 3.6 | 3.529411765 | 3.461538462 | 3.396226415 | 3.333333333 | 3.272727273 | 3.214285714 | 3.157894737 | 3.103448276 | 3.050847458 |
| 1/cos: | 1.001977173 | 1.001900275 | 1.00182778 | 1.001759358 | 1.001694709 | 1.00163356 | 1.001575664 | 1.001520793 | 1.00146874 | 1.001419316 |

| sides: | 60 | 61 | 62 | 63 | 64 | 65 | 66 | 67 | 68 | 69 |
|---|---|---|---|---|---|---|---|---|---|---|
| degrees: | 3 | 2.950819672 | 2.903225806 | 2.857142857 | 2.8125 | 2.769230769 | 2.727272727 | 2.686567164 | 2.647058824 | 2.608695652 |
| 1/cos: | 1.001372346 | 1.00132767 | 1.001285142 | 1.001244626 | 1.001205996 | 1.001169138 | 1.001133945 | 1.001100318 | 1.001068165 | 1.001037401 |

| sides: | 70 | 71 | 72 | 73 | 74 | 75 | 76 | 77 | 78 | 79 |
|---|---|---|---|---|---|---|---|---|---|---|
| degrees: | 2.571428571 | 2.535211268 | 2.5 | 2.465753425 | 2.432432432 | 2.4 | 2.368421053 | 2.337662338 | 2.307692308 | 2.278481013 |
| 1/cos: | 1.001007948 | 1.000979732 | 1.000952685 | 1.000926743 | 1.000901846 | 1.00087794 | 1.000854972 | 1.000832894 | 1.00081166 | 1.000791228 |

| sides: | 80 | 81 | 82 | 83 | 84 | 85 | 86 | 87 | 88 | 89 |
|---|---|---|---|---|---|---|---|---|---|---|
| degrees: | 2.25 | 2.222222222 | 2.195121951 | 2.168674699 | 2.142857143 | 2.117647059 | 2.093023256 | 2.068965517 | 2.045454545 | 2.02247191 |
| 1/cos: | 1.000771559 | 1.000752613 | 1.000734358 | 1.000716759 | 1.000699785 | 1.000683407 | 1.000667597 | 1.00065233 | 1.000637581 | 1.000623326 |

| sides: | 90 | 91 | 92 | 93 | 94 | 95 | 96 | 97 | 98 | 99 |
|---|---|---|---|---|---|---|---|---|---|---|
| degrees: | 2 | 1.978021978 | 1.956521739 | 1.935483871 | 1.914893617 | 1.894736842 | 1.875 | 1.855670103 | 1.836734694 | 1.818181818 |
| 1/cos: | 1.000609544 | 1.000596215 | 1.000583318 | 1.000570835 | 1.000558748 | 1.000547042 | 1.000535699 | 1.000524706 | 1.000514048 | 1.000503711 |

| sides: | 100 | 101 | 102 | 103 | 104 | 105 | 106 | 107 | 108 | 109 |
|---|---|---|---|---|---|---|---|---|---|---|
| degrees: | 1.8 | 1.782178218 | 1.764705882 | 1.747572816 | 1.730769231 | 1.714285714 | 1.698113208 | 1.682242991 | 1.666666667 | 1.651376147 |
| 1/cos: | 1.000493683 | 1.000483952 | 1.000474505 | 1.000465333 | 1.000456424 | 1.000447768 | 1.000439356 | 1.00043118 | 1.000423229 | 1.000415496 |